\documentclass{article}

\newtheorem{theorem}{Theorem}
\newtheorem{lemma}{Lemma}
\newtheorem{problem}{Problem}
\newcommand{\bt}{\begin{theorem}}
\newcommand{\et}{\end{theorem}}
\newcommand{\bl}{\begin{lemma}}
\newcommand{\el}{\end{lemma}}
\newcommand{\bp}{\begin{problem}}
\newcommand{\ep}{\end{problem}}
\newcommand{\pf}{{\bf Proof}.\ }
\newcommand{\eop}{\vspace{.8cm}}
\newcommand{\bq}{\begin{eqnarray*}}
\newcommand{\eq}{\end{eqnarray*}}
\newcommand{\be}{\begin{eqnarray}}
\newcommand{\ee}{\end{eqnarray}}
\newcommand{\beq}{\begin{equation}}
\newcommand{\eeq}{\end{equation}}
\newcommand{\benum}{\begin{enumerate}}
\newcommand{\eenum}{\end{enumerate}}
\newcommand{\ba}{\begin{array}}
\newcommand{\ea}{\end{array}}
\newcommand{\Z}{\ensuremath{\mathbf Z}}
\newcommand{\N}{\ensuremath{\mathbf N}}
\newcommand{\Q}{\ensuremath{\mathbf Q}}
\newcommand{\pol}{$\mathcal{F} = \{f_n(q)\}_{n=1}^{\infty}$}

\begin{document}
\title{Additive number theory and \\
the ring of quantum integers}
\author{Melvyn B. Nathanson\thanks{This work was supported
in part by grants from the NSA Mathematical Sciences Program
and the PSC-CUNY Research Award Program.}\\
Department of Mathematics\\
Lehman College (CUNY)\\
Bronx, New York 10468\\
Email: nathansn@alpha.lehman.cuny.edu}
\maketitle

\begin{abstract}
Let $m$ and $n$ be positive integers.
For the quantum integer $[n]_q = 1+q+q^2+\cdots + q^{n-1}$
there is a natural polynomial addition
such that $[m]_q \oplus_q [n]_q = [m+n]_q$
and a natural polynomial multiplication such that 
$[m]_q\otimes_q [n]_q = [mn]_q$.
These definitions lead to the construction
of the ring of quantum integers 
and the field of quantum rational numbers.
It is also shown that addition and multiplication of quantum integers 
are equivalent to elementary decompositions of intervals of integers
in additive number theory.
\end{abstract}

\section{Addition and multiplication}
Let $\N, \Z$ and \Q\ be the sets of positive integers, 
integers, and rational numbers, respectively.
We define the function
\[
[x]_q = \frac{1-q^x}{1-q}
\]
of two variables $x$ and $q$.
This is called the {\em quantum number} $[x]_q$.
Then
\[
[0]_q = 0,
\]
and for every positive integer $n$ we have
\[
[n]_q = \frac{1-q^n}{1-q} = 1 + q + \cdots + q^{n-1},
\]
which is the usual quantum integer $n$.
The negative quantum integers are
\bq
[-n]_q & = & \frac{1-q^{-n}}{1-q} 
= -\frac{1-q^n}{q^n(1-q)} \\
& = & -\frac{1}{q^n}[n]_q = -q^{-1}[n]_{q^{-1}}\\
& = & -\left( \frac{1}{q}+\frac{1}{q^2}+ \cdots +\frac{1}{q^n} \right).
\eq
Define quantum addition $\oplus_q$ as follows:
\[
[x]_q \oplus_q [y]_q = [x]_q + q^x[y]_q.
\]
Then 
\bq
[x]_q \oplus_q [y]_q & = & [x]_q + q^x[y]_q \\
& = & \frac{1-q^x}{1-q} + q^x\frac{1-q^y}{1-q} \\
& = & \frac{1-q^{x+y}}{1-q} \\
& = & [x+y]_q.
\eq
Define quantum multiplication $\otimes_q$ as follows:
\[
[x]_q \otimes_q [y]_q = [x]_q [y]_{q^x}
\]
Then
\bq
[x]_q \otimes_q [y]_q & = &[x]_q [y]_{q^x} \\
& = & \frac{1-q^x}{1-q} \frac{1-{q^x}^y}{1-q^x} \\
& = & \frac{1-q^{xy}}{1-q} \\
& = & [xy]_q.
\eq
The identities
\beq           \label{qiden}
[x]_q \oplus_q [y]_q = [x+y]_q \qquad\mbox{and}\qquad
[x]_q \otimes_q [y]_q = [xy]_q
\eeq
immediately imply that the set
\[
[\Z]_q = \{ [n]_q : n\in \Z\}
\]
is a commutative ring with the operations of quantum  addition $\oplus_q$
and quantum multiplication $\otimes_q.$  
The ring $[\Z]_q$ is called the {\em ring of quantum integers.}
The map $n \mapsto [n]_q$ from \Z\ to $[\Z]_q$  is a ring isomorphism.

For any rational number $m/n,$ the quantum rational number $[m/n]_q$ is
\bq
[m/n]_q & = & \frac{1-q^{m/n}}{1-q} 
= \frac{1-\left(q^{1/n}\right)^m}{1-\left(q^{1/n}\right)^n} \\
& = & \frac{  \frac{1-\left(q^{1/n}\right)^m}{1-q^{1/n}} } 
 {\frac{1-\left(q^{1/n}\right)^n}{1-q^{1/n}}} 
=  \frac{[m]_{q^{1/n}}} {[n]_{q^{1/n}}}.
\eq
Identities~(\ref{qiden}) imply that addition and multiplication 
of quantum rational numbers are well-defined.
We call
\[
[\Q]_q = \{ [m/n]_q : m/n \in \Q \}
\]
the {\em field of quantum rational numbers.}

If we consider $[x]_q$ as a function of real variables $x$ and $q$,
then
\[
\lim_{q\rightarrow 1} [x]_q = x
\]
for every real number $x$.

We can generalize the results in this section as follows:

\bt
Consider the function
\[
[x]_q = \frac{1-q^x}{1-q}
\]
in the variables $x$ and $q$.
For any ring $R$, not necessarily commutative, the set
\[
[R]_q = \{ [x]_q : x \in R\}
\]
is a ring with addition defined by 
\[
[x]_q \oplus_q [y]_q = [x]_q + q^x[y]_q.
\]
and multiplication by
\[
[x]_q \otimes_q [y]_q = [x]_q [y]_{q^x}
\]
The map from $R$ to $[R]_q$ defined by $x \mapsto [x]_q$ is a ring isomorphism.
\et

\pf
This is true for an arbitrary ring $R$ because the two identities 
in~(\ref{qiden}) are formal.

\section{Uniqueness of quantum arithmetic}
Let \pol\ be a sequence of polynomials in the variable $q$ 
that satisfies the addition and multiplication rules for quantum integers, that is, 
$\mathcal{F}$ satisfies the {\em additive functional equation}
\beq                    \label{afe}
f_{m+n}(q) = f_m(q) + q^mf_n(q)
\eeq
and the {\em multiplicative functional equation}
\beq                    \label{mfe}
f_{mn}(q) =f_m(q)f_n(q^m)
\eeq
for all positive integers $m$ and $n$.
Nathanson~\cite{nath02b} showed that there is a rich variety 
of sequences of polynomials that satisfy the multiplicative functional 
equation~(\ref{mfe}).
There is not yet a complete classification of solutions of~(\ref{mfe}), 
but there is a simple description of all solutions of the 
additive functional equation~(\ref{afe}).

\bt             \label{theorem:afe}
Let \pol\ be a sequence of functions that satisfies the
additive functional equation~(\ref{afe}).  Let $h(q) = f_1(q).$
Then
\beq               \label{afesol}
f_n(q) = h(q)[n]_q \qquad\mbox{for all $n \in \N.$}
\eeq
Conversely, for any function $h(q)$ the sequence of functions 
\pol\ defined by~(\ref{afesol}) is a solution of~(\ref{afe}).
In particular, if $h(q)$ is a polynomial in $q$, then $h(q)[n]_q$
is a polynomial in $q$ for all positive integers $n$, and all
polynomial solutions of~(\ref{afe}) are of this form.
\et

\pf
Suppose that \pol\ is a solution of the additive functional 
equation~(\ref{afe}).
Define $h(q) = f_1(q).$  Since $[1]_q = 1$ we have
\[
f_1(q) = h(q)[1]_q.
\]
Let $n \geq 2$ and suppose that $f_{n-1}(q) = h(q)[n-1]_q$.
From~(\ref{afe}) we have
\bq
f_n(q) & = & f_1(q)+qf_{n-1}(q) \\
& = & h(q)[1]_q + qh(q)[n-1]_q  \\
& = & h(q)([1]_q + q[n-1]_q)  \\
& = & h(q)[n]_q.
\eq
It follows by induction that $f_n(q) = h(q)[n]_q$ for all $n \in \N.$

Conversely, multiplying~(\ref{afe}) by $h(q)$, we obtain
\[
h(q)[m+n]_q = h(q)[m]_q + q^mh(q)[n]_q,
\]
and so the sequence $\{h(q)[n]_q\}_{n=1}^{\infty}$
is a solution of the additive functional equation~(\ref{afe}) for any function $h(q).$
This completes the proof.
\eop

We can now prove that the sequence of quantum integers is the only nontrivial solution
of the additive and multiplicative functional equations~(\ref{afe}) and~(\ref{mfe}).

\bt
Let \pol\ be a sequence of functions that satisfies both functional 
equations~(\ref{afe}) and~(\ref{mfe}).  Then either $f_n(q) = 0$ for all $n \in \N$
or $f_n(q) = [n]_q$ for all $n \in \N.$
\et

\pf
The multiplicative functional equation implies that
$f_1(q) = f_1(q)^2$, and so either $f_1(q) = 0$ or $f_1(q) = 1.$
Since \pol\ also satisfies the additive functional equation,
it follows from Theorem~\ref{theorem:afe} that either $f_n(q) = 0$ for all $n$
or $f_n(q) = [n]_q$ for all $n$.
This completes the proof.
\eop

\section{Additive number theory}
In this section we show that the addition and multiplication rules
for quantum integers correspond to elementary decompositions of 
finite sets of integers in additive number theory.

Let $A$ and $B$ be sets of integers, and let $m$ be an integer.  
We define the {\em dilation}
\[
m\ast A = \{ma : a\in A \},
\]
the {\em translation}
\[
m + A = \{m+a : a\in A \},
\]
and the {\em sumset}
\[
A+B = \{a+b : a\in A, b\in B\}.
\]
We write
\[
A\oplus B = C
\]
if $A+B = C$ and every integer in $C$ can be written 
uniquely in the form $a+b$ for some $a \in A$ and $b\in B.$

In additive number theory we consider partitions of a set of integers 
into a disjoint union of subsets, and decompositions of a set of integers 
into a sum of sets of integers.  
Denote by $[n]$ the set of the first $n$ nonnegative integers, that is,
\[
[n] = \{0,1,2,\ldots, n-1\}.
\]
We have the partition
\beq              \label{m+n}
[m+n] = [m] \cup (m + [n]),\qquad\mbox{where $[m] \cap (m+[n]) = \emptyset$,}
\eeq
and the direct sum decomposition
\beq          \label{mn}
[mn] = [m] \oplus m\ast[n].
\eeq
If $m_1,\ldots, m_r$ are positive integers, then, by induction, 
we have the partition
\beq              \label{m1+mr}
[m_1 + m_2 + \cdots + m_r] = \bigcup_{j=1}^r \left( \sum_{i=1}^{j-1}m_i + [m_j] \right)
\eeq
into pairwise disjoint sets, and the direct sum decomposition
\beq          \label{m1mr}
[m_1m_2\cdots m_r] = \bigoplus_{j=1}^r \left( \prod_{i=1}^{j-1}m_i \ast [m_j] \right).
\eeq

To each finite set $A$ of integers we associate the Laurent polynomial
\[
F_A(q) = \sum_{a\in A} q^a.
\]
This is called the {\em generating function} for $A$.
From the definitions of dilation, translation, and sumset,
we have the generating function identities
\[
F_{m\ast A}(q) = F_A(q^m),
\]
\[
F_{m+A}(q) = q^mF_A(q),
\]
and
\[
F_{A\oplus B}(q) = F_A(q)F_B(q).
\]
If $A \cap B = \emptyset,$ then
\[
F_{A\cup B}(q) = F_A(q) + F_B(q).
\]

The generating function for the set $[n]$ is the quantum integer $[n]_q$, since
\[
F_{[n]}(q) = 1 + q + \cdots + q^{n-1} = [n]_q.
\]
Rewriting the partition identity~(\ref{m+n}) in terms of generating functions, 
we obtain
\bq
[m+n]_q 
& = & F_{[m+n]}(q) \\
& = & F_{[m] \cup (m + [n])}(q) \\
& = & F_{[m]}(q)+ F_{m+[n]}(q)\\
& = & F_{[m]}(q)+ q^mF_{[n]}(q)\\
& = & [m]_q + q^m[n]_q.
\eq
The sumset decomposition~(\ref{mn}) of the interval $[mn]$ gives
\bq
[mn]_q 
& = & F_{[mn]}(q) \\
& = & F_{[m] \oplus m\ast[n]}(q) \\
& = & F_{[m]}(q)F_{m\ast [n]}(q) \\
& = & F_{[m]}(q)F_{[n]}(q^m) \\
& = & [m]_q [n]_{q^m}.
\eq
Similarly, the additive number theoretic identities~(\ref{m1+mr}) 
and~(\ref{m1mr}) yield the quantum integer identities
\[
[m_1 + m_2 + \cdots + m_r]_q = \sum_{j=1}^r q^{\sum_{i=1}^{j-1}m_i} [m_j]_q 
\]
and
\[
[m_1m_2\cdots m_r]_q = \prod_{j=1}^r  [m_j]_{q^{\prod_{i=1}^{j-1}m_i}}.
\]
In this way we see that the addition and multiplication rules 
for quantum integers are equivalent to elementary statements in additive number theory.

\end{document}